\documentclass[twoside,11pt]{amsart}
\usepackage{amsmath,latexsym,amssymb,times,mathptm,amsfonts,enumerate}

\def\proof{\noindent{\bf{Proof.} }}
\def\sqr#1#2{{\vcenter{\hrule height.#2pt
        \hbox{\vrule width.#2pt height#1pt \kern#1pt
                \vrule width.#2pt}
        \hrule height.#2pt}}}
\def\square{\mathchoice\sqr64\sqr64\sqr{4}3\sqr{3}3}
\def\QED{\hfill$\square$}

\def\tratto{\mbox{\rule{2mm}{.2mm}$\;\!$}}

\newtheorem{theorem}{Theorem}[section]

\newtheorem{lemma}[theorem]{Lemma}
\newtheorem{proposition}[theorem]{Proposition}

\newtheorem{conjecture}[theorem]{Conjecture}

\newtheorem{definition}[theorem]{Definition}
\newtheorem{remark}[theorem]{Remark}

\newtheorem{example}[theorem]{Example}
\newtheorem{question}[theorem]{Question}

\newcommand{\s}{\bigskip}

\newcommand{\f}[1]{\ensuremath{\mathfrak{#1}}}
\newcommand{\ol}[1]{\ensuremath{\overline{#1}}}

\newcommand{\height}{\ensuremath{{\rm{ht}}\;}}
\newcommand{\depth}{\ensuremath{{\rm{depth}}\;}}
\newcommand{\core}[1]{\ensuremath{{\rm{core}}(#1)}}
\newcommand{\ul}[1]{\ensuremath{\underline{#1}}}
\numberwithin{equation}{section}

\newcommand{\ds}{\displaystyle}

\begin{document}
\baselineskip=16pt

\title{Computing the core of ideals in arbitrary characteristic}
\author[L. Fouli]{Louiza Fouli}
\address{Department of Mathematics, 1 University Station C1200, University of Texas, Austin, Texas 78712, USA}
\email{lfouli@math.utexas.edu}

\vspace{-0.1in}

\begin{abstract}
Let $R$ be a local Gorenstein ring with infinite residue field of
arbitrary characteristic. Let $I$ be an $R$--ideal with $g=\height I
>0$, analytic
  spread $\ell$, and let $J$ be a minimal reduction of $I$. We further assume that $I$ satisfies $G_{\ell}$ and ${\depth}
  \, R/I^j \geq \dim R/I -j+1$ for $1 \leq j \leq \ell-g$. The question we are interested in is whether $\core{I}=J^{n+1}:\ds \sum_{b \in I} (J,b)^n$ for $n \gg 0$. In
the case of analytic spread one Polini and Ulrich show that this is
true with even weaker assumptions (\cite[Theorem 3.4]{PU}). We give
a negative answer to this question for higher analytic spreads and
suggest a formula for the core of such ideals.
\end{abstract}

\vspace{-0.2in}

\maketitle

\section{Introduction}
Throughout let $R$ be a Noetherian ring. If $R$ is a Noetherian
local ring with maximal ideal $\f{m}$ then we denote the residue
field of $R$ by $k=R/\f{m}$. Let $I$ be an $R$--ideal. In order to
study an ideal $I$, Northcott and Rees introduced the notion of a
reduction of an ideal. A reduction in general is a simplification of
the ideal itself. Recall that a {\it reduction} of an ideal $I$ is a
subideal $J$ such that $I^{n+1}=JI^{n}$, for some nonnegative
integer $n$ (\cite{NR}). This condition is equivalent to $I$ being
integral over $J$. Moreover, reductions preserve a number of
properties of the ideal and thus it is customary to shift the
attention from the ideal to its reductions. In the case that $R$ is
a Noetherian local ring we may consider {\it minimal reductions},
which are minimal with respect to inclusion. Northcott and Rees
prove that if the residue field $k$ of $R$ is infinite then
minimal reductions do indeed exist and they correspond to Noether
normalizations of the special fiber ring
$\mathcal{F}(I):={\underset{i \geq 0}\oplus}
I^{i}/\f{m}I^{i}=R/\f{m}\oplus I/\f{m}I\oplus \ldots \oplus
I^{i}/\f{m}I^{i} \oplus \ldots$ of $I$ (\cite{NR}). In particular
this shows that minimal reductions are not unique. Recall that the
{\it analytic spread of } $I$, $\ell(I)$, is the Krull dimension of
the special fiber ring $\mathcal{F}(I)$, i.e., $\ell=\ell(I)=\dim
\mathcal{F}(I)$. If $k$ is infinite Northcott and Rees also show that
for any minimal reduction $J$ of $I$ one has $\mu(J)=\ell(I)$, where
$\mu(J)$ denotes the minimal number of generators of $J$
(\cite{NR}).

In order to counteract the lack of uniqueness of minimal reductions
Rees and Sally consider the intersection over all (minimal)
reductions, namely the {\it core} of the ideal (\cite{RS}). Then
$\core I={\ds {\ds \bigcap_{J}^{} }  J}$, where $J$ is a (minimal)
reduction of $I$. The core arises naturally in the context of
Brian\c{c}on--Skoda kind of theorems. If $R$ is a regular local ring
of dimension $d$ and $I$ is an $R$--ideal, then the
Brian\c{c}on--Skoda theorem states that $\overline{I^{d}} \subset
J$, for every reduction $J$ of $I$, or equivalently
$\overline{I^{\;d}} \subset \core I$, where `${}^{\tratto}$' denotes
the integral closure of the corresponding ideal. Huneke and Swanson
(\cite{HS}) showed a connection between the work of Lipman
(\cite{L}) on the adjoint of an ideal and the core. The core is a
priori an infinite intersection. Hence there is significant
difficulty in computing this ideal. The question of finding explicit
formulas that compute the core has been addressed in the work of
Corso, Huneke, Hyry, Polini, Smith, Swanson, Trung, Ulrich and
Vitulli (\cite{CPU01}, \cite{CPU02}, \cite{HS}, \cite{HT},
\cite{HySm1}, \cite{HySm2}, \cite{PU}, \cite{PUV}). Moreover, Hyry
and Smith have discovered a connection with a conjecture by Kawamata
on the non--vanishing of sections of line bundles (\cite{HySm2}).

In this paper we are primarily interested in a formula for the core
of an ideal shown by Polini and Ulrich which states:
\begin{theorem}{\rm{{(\cite[Theorem~4.5]{PU})}}}\label{formulaPU}
Let $R$ be a local Gorenstein ring with infinite residue field $k$,
let $I$ be an $R$--ideal with $g= \height I > 0$ and $\ell=
\ell(I)$, and let $J$ be a minimal reduction of $I$ with reduction
number $r$. Assume $I$ satisfies $G_{\ell}$ and $\depth R/I^j \geq
{\rm dim} \, R/I -j+1$ for $1 \leq j \leq \ell-g$, and either ${\rm
char} \, k =0 $ or ${\rm char} \, k  > r -\ell + g $. Then
\[
{\rm core}(I)= J^{\;n+1} : I^n
\]
for every $n \geq {\rm max} \{\; r - \ell + g, 0 \;\}$.
\end{theorem}

The goal will be clear once it is understood how the
formula in Theorem~\ref{formulaPU} arises. In general Polini and
Ulrich show that:

\begin{theorem}\label{geninclu}{\rm {(\cite[Remark~4.8]{PU})}
Let $R$ be a local Gorenstein ring with infinite residue field, let
$I$ be an $R$--ideal with $g=\height I>0$ and $\ell=\ell(I)$, and
let $J$ be a minimal reduction of $I$ with reduction number $r$. Assume $I$ satisfies
$G_{\ell}$ and $\depth R/I^j \geq \dim R/I -j+1$ for $1 \leq j \leq
\ell-g$. Then
\begin{equation}\label{general inclusion}
J^{n+1}:I^{n} \subset {\rm{core}}(I) \subset J^{n+1}:\ds \sum_{b \in
I} (J,b)^n
\end{equation}
for every $n \geq {\rm{max}} \{\;r-\ell+g,0\;\}$.}
\end{theorem}

These inclusions hold in any characteristic. The condition on the
characteristic of the residue field in Theorem~\ref{formulaPU}
implies that the two bounds for the core in equation~(\ref{general
inclusion}) coincide. This gives the formula in
Theorem~\ref{formulaPU}.

When the analytic spread of $I$ is one, Polini and Ulrich also show
the following:

\begin{theorem}{\rm{{(\cite[Theorem~3.4]{PU})}}}\label{ell=one}
Let $R$ be a local Cohen--Macaulay ring with infinite residue field,
let $I$ be an $R$--ideal with $\ell(I)= \height I = 1$, and let $J$ be a minimal reduction of $I$. Then for $n \gg 0$
\[
 {\rm core} (I) =  J^{n+1} :
\sum_{b \in I} (J, b)^n.
\]
\end{theorem}

Notice that Theorem~\ref{ell=one} holds in any characteristic. In
the same paper Polini and Ulrich also exhibit a class of examples
where $\ell(I)=1$ and ${\rm{core}}(I) \neq J^{\;n+1}:I^{n}$
({\cite[Example~4.9]{PU}}). Thus a natural question arises:

\begin{question}\label{core=Ln?}{\rm
Under the same assumptions as in Theorem \ref{formulaPU}, except for
the condition on the characteristic of $k$, is ${\rm{core}}(I)
=J^{n+1}:{\ds \sum_{b \in I} (J,b)^n}$ for some $n \gg 0$?}
\end{question}

The purpose of this paper is to answer Question~\ref{core=Ln?}. In
order to answer this question we first seek to better understand the
ideal $J^{n+1}:{\ds \sum_{b \in I} (J,b)^n}$. One of the
difficulties lies in computing ${\ds \sum_{b \in I} (J,b)^n}$. We
devote Section~\ref{The ideal $K_n$} to understanding this ideal. In
Theorem~\ref{Kn general} we give an explicit algorithm for computing
this ideal. Once we are able to compute it we are interested in the
behaviour of the ideal $J^{n+1}:{\ds \sum_{b \in I} (J,b)^n}$, which
we address in Section~\ref{The ideal $L_n$}. In
Section~\ref{Examples} we finally answer Question~\ref{core=Ln?}.

Before we proceed any further we need to explain some of the
conditions that are used in Theorem~\ref{formulaPU} and throughout
this paper. Let $R$ be a Noetherian ring, $I$ an $R$--ideal and $s$
an integer. We say that $I$ satisfies $G_s$ if $\mu(I_{\f{p}}) \leq
\dim R_{\f{p}}$ for every $\f{p} \in V(I)$ with $\dim R_{\f{p}} \leq
s-1$. If $R$ is a Noetherian local ring of dimension $d$ and $\f{m}$
is the maximal ideal of $R$, then any $\f{m}$--primary ideal
satisfies $G_{d}$. An additional technical condition that is
connected with the study of the core is the assumption ${\rm{depth}}
\, R/I^{j} \geq {\rm{dim}} \, R/I-j+1$ for $1 \leq j \leq \ell-g$,
where $g$ is the height of $I$ and $\ell$ is the analytic
spread of $I$.

Let $R$ be a local Gorenstein ring with maximal ideal $\f{m}$ and
infinite residue field and let $I$ be an $R$--ideal with height $g$
and analytic spread $\ell$. Then $I$ satisfies $G_{\ell}$ and ${\rm
depth} \, R/I^j \geq \dim R/I -j+1$ for $1 \leq j \leq \ell -g$ in
the following cases:
\begin{enumerate}
\item[{\rm (a)}] $I$ is $\f{m}$--primary, or more generally $I$ is equimultiple which means $\ell=g$.

\item[{\rm (b)}] $I$ is a one--dimensional generic complete intersection
ideal, or more generally $I$ is a generic complete intersection
Cohen--Macaulay ideal with $\ell \leq g+1$ (\cite[p. 259]{AvHer}.
\end{enumerate}
In the presence of the $G_{\ell}$ property the depth condition on
the powers of $I$ as above is satisfied if $I$ is perfect of height
$2$, or if $I$ is perfect Gorenstein of height $3$, or more
generally if $I$ is in the linkage class of a complete intersection
ideal (licci) (\cite[1.11]{Hun}).

An interesting invariant of an ideal $I$ is the reduction number of
$I$. The {\it reduction number of} $I$ {\it with respect to} $J$ is
the integer $r_{J}(I)={\rm min} \{\; n \;|\; I^{n+1}=JI^{n}\;\}$,
where $J$ is a reduction of $I$. The {\it reduction number of} $I$,
$r(I)$, is defined to be ${\rm min}\; \{ \; r_{J}(I)\;|\; J {\rm \;
minimal \; reduction\;  of}\; I\;\}$. The reduction number of $I$ is
connected with the study of blowup algebras and their
Cohen-Macaulayness.

\section{The ideal $K_n$}{\label{The ideal $K_n$}}

Let $R$ be a Noetherian local ring with infinite residue field and
let $I$ be an $R$--ideal. Let $J$ be a minimal reduction of $I$. As
a starting point of our work we first seek to better understand the
ideal ${\ds \sum_{b \in I} (J,b)^n }$ as it is connected with the
core of $I$ by work of Polini and Ulrich (\cite{PU}). Our first goal
is to find an efficient way to compute this ideal. We start
investigating such an ideal in a general setting.

\begin{definition}
{\rm Let $R$ be a Noetherian ring and let $J \subset I$ be
$R$--ideals. Let $n$ be a positive integer. We denote by
$K_{n}(J,I)$ the $R$--ideal $\ds \sum_{b \in I} (J,b)^n$ and by
$L_{n}(J,I)$ the $R$--ideal $J^{n+1}:K_{n}(J,I)$. When the ideal $I$
is understood we will denote these ideals by $K_{n}(J)$ and
$L_{n}(J)$, respectively. If in addition the ideal $J$ is understood
then we will use $K_{n}$ and $L_{n}$, respectively.}
\end{definition}

The following lemma gives an explicit description of a (not
necessarily minimal) generating set for $K_n(J,I)$.

\begin{lemma} \label{lemmagens}
Let $R$ be a Noetherian local ring with infinite residue field. Let
$T$ be an algebra over $R$ and let $J \subset I$ be $T$--ideals.
Assume that $J=(f_{1},\ldots , f_t)$ and $I=(J, f_{t +1}, \ldots
,f_{m})$. Then
\[
K_n(J,I)=(\{{\ds \prod_{j=t+1}^{m-1} {{n-{
\overset{j-1}{\underset{i=1} \sum}
\nu_i}}\choose{\nu_j}}f_{1}^{\nu_{1}}\cdots f_{m}^{\nu_{m}} } \} ),
\]
where $\nu_{1}, \nu_{2}, \cdots ,\nu_{m}$ range over all
nonnegative integers with $\nu_{1}+ \nu_{2} + \cdots +\nu_{m}=n$.
\end{lemma}

\proof Let $k$ denote the residue field of $R$. Let $A$ be the ideal
$(\{{\ds \prod_{j=t+1}^{m-1} {{n-{\overset{j-1}{\underset{i=1} \sum}
\nu_i}}\choose{\nu_j}}f_{1}^{\nu_{1}}\cdots f_{m}^{\nu_{m}} } \} )$,
where $\nu_{1}+ \nu_{2} + \cdots +\nu_{m}=n$.

We first show that $K_n=K_{n}(J,I) \subset A$. It suffices to check
that the generators of $K_{n}(J,I)$ are in $A$. Let $f$ be such a
generator. We may assume that $f$ is of the form $f_{1}^{\nu_{1}}
\cdots f_{t}^{\nu_{t}}y^{n-s}$, where $\nu_{1}+ \cdots +\nu_{t}=s$,
$0 \leq s \leq n$, and $y=\overset{m}{\underset{i=t+1}
\sum}f_{i}g_{i}$, with $g_{i} \in T$ for $t+1 \leq i \leq m$. Then
\[\begin{array}{ll}
f&=f_{1}^{\nu_{1}} \cdots f_{t}^{\nu_{t}}(f_{t+1}g_{t+1}+ \cdots +f_{m}g_{m})^{n-s}\\
&= {\underset{\ul{\nu}}\sum}\beta_{\ul{\nu}}f_{1}^{\nu_{1}} \cdots
f_{t}^{\nu_{t}}f_{t+1}^{\nu_{t+1}} \cdots
f_{m}^{\nu_{m}}g_{t+1}^{\nu_{t+1}} \cdots g_{m}^{\nu_{m}},
\end{array}\]
where
\begin{enumerate}
\item[{\rm (a)}] $\ul{\nu}=(\nu_{t+1}, \cdots ,\nu_{m})$;
\item[{\rm (b)}] $\nu_{t+1}+ \cdots +\nu_{m}=n-s$ and $\nu_{1}+ \cdots
+\nu_{t}=s$, \quad and
\item[{\rm (c)}] $\beta_{\ul{\nu}}={\ds \prod_{j=t+1}^{m-1}} {\ds {{n- {\overset{j-1}{\underset{i=1} \sum}\nu_{i}}}
\choose{\nu_{j}}}}$.
\end{enumerate}
Hence  clearly $f \in A$ and $K_n \subset A$. It remains to show
that $A \subset K_{n}$.

Fix $\nu_{1}, \cdots, \nu_{t}$. Let $N=n-{\ds \sum_{i=1}^{t}}
\nu_{i}$. As above write $\beta_{\ul{\nu}}={\ds \prod_{j=t+1}^{m-1}}
{\ds {{n- {\overset{j-1}{\underset{i=1}
\sum}\nu_{i}}}\choose{\nu_{j}}}}$, where
$\underline{\nu}=(\nu_{t+1},\cdots, \nu_{m})$ and $\nu_{t+1}+
\cdots+ \nu_{m} =N$. We prove that
\begin{equation*}\label{(1)}
\beta_{\underline{\nu}}f_{1}^{\nu_{1}}\cdots f_{m}^{\nu_{m}} \in
K_{n}(J).
\end{equation*}

Let $M_{1}, \cdots, M_{D}$ be the monomials of degree $N$ in the
variables $X_{t+1}, \cdots, X_{m}$. Then the elements
$\beta_{\underline{\nu}}f_{1}^{\nu_{1}}\cdots f_{m}^{\nu_{m}}$
become $\beta_{i}f_{1}^{\nu_{1}} \cdots
f_{t}^{\nu_{t}}M_{i}(\ul{f})$, where $\underline{f}=(f_{t+1}, \cdots
,f_{m})$ and $1 \leq i \leq D$. We choose $\alpha_{t+1}, \cdots,
\alpha_{m}$ elements in $R$ such that for all $1\leq i \leq D$ the
images of $M_{i}(\ul{\alpha})$ in $k$ are distinct, where
$\ul{\alpha}=(\alpha_{t+1}, \cdots, \alpha_{m})$. Since the residue
field $k$ is infinite it is possible to choose such elements. Then
for $ 0 \leq e \leq D-1$ we have
\begin{equation}\label{(2)}
K_{n} \ni f_{1}^{\nu_{1}} \cdots f_{t}^{\nu_{t}}\left({\ds
\sum_{j=t+1}^{m}}\alpha_{j}^{e}f_{j}\right)^{N} =
f_{1}^{\nu_{1}}\cdots f_{t}^{\nu_{t}}{\ds \sum_{i=1}^D}
\beta_{i}(M_{i}(\underline{\alpha}))^{e}M_{i}(\underline{f}).
\end{equation}

Let $B$ denote the $D \times D$ matrix whose $(j,i)$ entry is
$(M_{i}(\underline{\alpha}))^{j-1}$ and let $C$ denote the $D \times
1$ matrix whose $i^{\rm{th}}$ entry $\beta_{i}f_{1}^{\nu_{1}} \cdots
f_{t}^{\nu_{t}}$ is $M_{i}(\underline{f})$. The entries of $BC$ are
in $K_{n}$ according to equation~(\ref{(2)}). Notice that $B$ is a
Vandermondt matrix and hence the determinant of it is the product of
the differences of all $M_{i}(\underline{\alpha})$. By the choice of
$\underline{\alpha}$ these differences have non--zero images in $k$,
and therefore are units in $R$ and thus in $T$. This implies that
${\rm{det}} \,B$ is a unit and hence $B$ is invertible. Therefore
the entries of $C$ are in $K_{n}$. \QED

\s

\begin{definition}\label{Gen elts}{\rm
Let $R$ be a Noetherian local ring with infinite residue field $k$.
Let $I=(f_{1},\cdots ,f_{m})$ be an $R$--ideal and let $t$ be a
fixed positive integer. We say that $b_{1}, \cdots, b_{t}$ are
$t$ \textit{general elements in} $I$ if there exists a dense open
subset $U$ of $\mathbb{A}_{k}^{tm}$ such that for $1 \leq i \leq t$
and $1 \leq j \leq m$ we have that $b_{i}={\ds \sum_{j=1}^m}
\lambda_{ij}f_{j}$, where
$\underline{\underline{\lambda}}=[\lambda_{ij}]_{ij} \in
\mathbb{A}_{R}^{tm}$ and $\underline{\overline{\underline{\lambda}}}
\in U$ vary in $U$, where $\ul{\ol{\ul{\lambda}}}$ is the image of
$\ul{\ul{\lambda}}$ in $\mathbb{A}_{k}^{tm}$.

The ideal $J$ is called a \textit{general minimal reduction} of $I$
if $J$ is a reduction of $I$ generated by $\ell(I)$ general elements
in $I$.}
\end{definition}

\s

\begin{theorem} \label{Kn general}
Let $R$ be a Noetherian local ring with infinite residue field. Let
$J \subset I$ be $R$--ideals. Let $n$ be a positive integer and
$K_{n}=\ds \sum_{b \in I} (J,b)^n$. For a fixed integer $t$ let
$b_{1}, \cdots, b_{t+1}$ be $t+1$ general elements in $I$. Set
$C_j={\ds \sum_{i=1}^j }(J,b_{i})^n$ for $1\leq j \leq t+1$. Assume
that $C_{t}=C_{t+1}$. Then\[K_n=\sum_{i=1}^t (J,b_{i})^n= C_{t}.\]
\end{theorem}

\proof Assume that $J=(f_1,
\cdots ,f_{s})$ and $I=(f_1, \cdots , f_{s}\;, \;f_{s+1}, \cdots,
f_{m})$. Let $k$ denote the residue field of $R$. Clearly $C_{t}={\ds \sum_{i=1}^t}(J,b_{i})^n \subset
K_{n}$. Notice that there exists a positive integer $t ' >t$ such
that $K_{n}=\ds \sum_{b \in I}(J,b)^n={\ds \sum_{i=1}^{t'}
}(J,y_{i})^n $ for some $y_{i} \in I$.

We consider the natural projection maps $\pi_{j}:
\mathbb{A}_{k}^{t'm} \rightarrow \mathbb{A}_{k}^{(t+1)m}$ where:
\begin{enumerate}
\item[{\rm (a)}] $t+1 \leq j \leq t'$;

\item[{\rm (b)}] $\pi_{j}((\ul{a_{1}}, \cdots, \ul{a_{t}}, \ul{a_{t+1}},
\cdots, \ul{a_{j}}, \cdots, \ul{a_{t'}}))=(\ul{a_{1}}, \cdots,
\ul{a_{t}}, \ul{a_{j}})$, and

\item[{\rm (c)}] $\ul{a_{i}} \in \mathbb{A}_{k}^{m}$.
\end{enumerate}

Let $\underline{\underline{\lambda}}= [\lambda_{ij}]_{i,j} \in
\mathbb{A}_{R}^{(t+1)m}$ and let
$\ol{\underline{\underline{\lambda}}}$ denote the image of
$\underline{\underline{\lambda}}$ in $\mathbb{A}_{k}^{(t+1)m}$. By
our assumption there exists a dense open subset $U \subset
\mathbb{A}_{k}^{(t+1)m}$ such that $C_{t+1}={\ds
\sum_{i=1}^{t+1}}(J,b_{i})^n={\ds \sum_{i=1}^{t}}(J,b_{i})^n=C_{t}$,
where $b_{i}={\ds \sum_{j=1}^m} \lambda_{ij} f_{j}$ for $1 \leq i
\leq t+1$ and $\ol{\underline{\underline{\lambda}}} \in U$. Notice
that $V={\ds \bigcap_{j=t+1}^{t'} \pi_{j}^{-1}(U)}$ is a dense open
subset in $\mathbb{A}_{k}^{t'm}$. For $1 \leq i \leq t'$ let
$b_{i}={\ds \sum_{j=1}^m} \lambda_{ij} f_{j}$, where
$\ol{\underline{\underline{\lambda}}} \in V$. By the construction of
$V$ one has that $C_{t'}=C_{t}$. So it suffices to show that
$K_{n}=C_{t'}$.

Let $T=R[X_{ij}]$, where $1\leq i \leq t'$ and $1 \leq j \leq m$.
For $1\leq i \leq t'$ let $y_{i}={\ds \sum_{j=1}^m} \lambda_{ij}^{0}
f_{j}$ with $\lambda_{ij}^{0} \in \mathbb{R}$, and write
$\underline{\underline{\lambda}}^0= [\lambda_{ij}^{0}]_{i,j}$ and
$\underline{\underline{X}}=[X_{ij}]_{i,j}$. Consider the $T$--ideal
$\widetilde{K_{n}}={\ds \sum_{i=1}^{t'}}(JT,{\ds
\sum_{j=1}^{m}}X_{ij}f_{j})^n$ and the $R$--homomorphisms $\pi_{
\underline{\underline{\lambda}}} : T \rightarrow R$ that send
$\underline{\underline{X}}$ to $ \underline{\underline{\lambda}}$,
where $\underline{\underline{\lambda}} \in \mathbb{A}_{R}^{t'm}$.
Notice that $\pi_{ \underline{\underline{\lambda}}^0}(
\widetilde{K_{n}})=K_{n}$. Therefore we have
\begin{equation}\label{(3)}
K_{n}T \subset
\widetilde{K_{n}}+(\underline{\underline{X}}-\underline{\underline{\lambda}}^0),
\end{equation}
and
\begin{equation}\label{(4)}
\widetilde{K_{n}}\subset K_{n}(JT,IT) = K_{n}T,
\end{equation}
where the last equality holds by Lemma \ref{lemmagens}.

Write $\f{m}_{\underline{\underline{\lambda}}}$ for the maximal
ideals
$(\f{m},\underline{\underline{X}}-\underline{\underline{\lambda}})$
of $T$. Localizing equation (\ref{(4)}) at these maximal ideals
gives
\begin{equation*}\label{(5)}
(\widetilde{K_{n}})_{\f{m}_{\underline{\underline{\lambda}}}}
\subset K_{n}T_{\f{m}_{\underline{\underline{\lambda}}}},
\end{equation*}
and combining this with equation (\ref{(3)}) yields
\begin{equation}\label{(6)}
K_{n}T_{\f{m}_{\underline{\underline{\lambda}}^{0}}} =
(\widetilde{K_{n}})_{\f{m}_{\underline{\underline{\lambda}}^{0}}}+
(\underline{\underline{X}}-\underline{\underline{\lambda}}^0) \cap
K_{n}T_{\f{m}_{\underline{\underline{\lambda}}^{0}}}.
\end{equation}
Since
$\underline{\underline{X}}-\underline{\underline{\lambda}}^{0}$ is a
regular sequence on
$\ds{T_{\f{m}_{\underline{\underline{\lambda}}^{0}}}/{K_{n}T_{\f{m}_{\underline{\underline{\lambda}}^{0}}}}}$,
equation (\ref{(6)}) becomes
\begin{equation*}\label{(7)}
K_{n}T_{\f{m}_{\underline{\underline{\lambda}}^{0}}} =
(\widetilde{K_{n}})_{\f{m}_{\underline{\underline{\lambda}}^{0}}}+
(\underline{\underline{X}}-\underline{\underline{\lambda}}^0)
K_{n}T_{\f{m}_{\underline{\underline{\lambda}}^{0}}}
\end{equation*} and thus
\begin{equation}\label{(8)}
K_{n}T_{\f{m}_{\underline{\underline{\lambda}}^{0}}}=(\widetilde{K_{n}})_{\f{m}_{\underline{\underline{\lambda}}^{0}}}
\end{equation} by Nakayama's lemma.

Equation (\ref{(8)}) allows us to conclude that
$M_{\f{m}_{\underline{\underline{\lambda}}^{0}}}=0$ for the
$T$--module
$M=\ds{\frac{\widetilde{K_{n}}+K_{n}T}{\widetilde{K_{n}}}}$. Hence
$\f{m}_{\underline{\underline{\lambda}}^{0}} \notin {\rm{Supp}}(M)$
and thus there exists a dense open set $U \subset
\mathbb{A}_{k}^{mt'}$ such that
$M_{\f{m}_{\underline{\underline{\lambda}}}}=0$ for all
$\ol{\underline{\underline{\lambda}}} \in U$, where
$\underline{\underline{\lambda}} \in \mathbb{A}_{R}^{mt'}$ and
$\ol{\underline{\underline{\lambda}}}$ denotes the image of
$\underline{\underline{\lambda}}$ in $\mathbb{A}_{k}^{mt'}$.
Therefore
$\ds{\frac{(\widetilde{K_{n}})_{\f{m}_{\underline{\underline{\lambda}}}}+K_{n}T_{\f{m}_{\underline{\underline{\lambda}}}}}
{(\widetilde{K_{n}})_{\f{m}_{\underline{\underline{\lambda}}}}}}=0$
for all $\ol{\underline{\underline{\lambda}}} \in U$. In other
words, for all $\ol{\underline{\underline{\lambda}}} \in U$ we have
%\begin{equation*}\label{(9)}
$(\widetilde{K_{n}})_{\f{m}_{\underline{\underline{\lambda}}}}+K_{n}T_{\f{m}_{\underline{\underline{\lambda}}}}=
(\widetilde{K_{n}})_{\f{m}_{\underline{\underline{\lambda}}}}$.
%\end{equation*}
Let
$\rho_{\underline{\underline{\lambda}}}:T_{\f{m}_{\underline{\underline{\lambda}}}}
\rightarrow R$ be the $R$--homomorphism that sends
$\underline{\underline{X}}$ to $\underline{\underline{\lambda}}$.
Then the above equation yields $\rho_{
\underline{\underline{\lambda}}}((\widetilde{K_{n}})_{\f{m}_{\underline{
 \underline{_\lambda}}}}+K_{n}T_{\f{m}_{\underline{\underline{\lambda}}}})=
\rho_{\underline{\underline{\lambda}}}((\widetilde{K_{n}})_{\f{m}_{\underline{\underline{\lambda}}}})$
for all $\ol{\underline{\underline{\lambda}}} \in U$. Hence
$\pi_{\ul{\ul{\lambda}}}(\widetilde{K_{n}})+K_{n}=\pi_{\ul{\ul{\lambda}}}(\widetilde{K_{n}})$
for all $\ol{\underline{\underline{\lambda}}} \in U$. As
$\pi_{\ul{\ul{\lambda}}}(\widetilde{K_{n}})\subset K_{n}$ we
conclude that
$K_{n}=\pi_{\ul{\ul{\lambda}}}(\widetilde{K_{n}})=C_{t'}$ for all
$\ol{\underline{\underline{\lambda}}}$ in a suitable dense open
subset of $\mathbb{A}_{k}^{t'm}$. \QED

\s

\begin{remark}{\rm
Notice that Theorem~\ref{Kn general} provides an algorithm for
computing the ideal $K_{n}$ for any positive integer $n$. We apply
this algorithm in computations using the computer algebra program
Macaulay 2 (\cite{M2}) in Section~\ref{Examples}.}
\end{remark}

\s

\section{The ideal $L_n$}{\label{The ideal $L_n$}}
In light of the algorithm given in Theorem~\ref{Kn general} we are
now able to compute the ideals $L_n=L_n(J)=J^{n+1}:\ds{ \sum_{b \in
I} (J,b)^n}$. Recall that our goal is to determine whether
$\core{I}=L_n$ for $n\gg 0$ (Question~\ref{core=Ln?}). However,
determining what $n \gg 0$ means is another challenge of its own. If
$\core{I} \neq L_n$ for some $n$ then in principle there is still a
possibility that $\core{I}=L_m$ for $m>n$. We need to determine how
one can effectively decide when $\core{I} \neq L_{n}$ for all $n
>0$. In this section we will prove that the ideals $L_{n}$ stabilize
past a computable integer (Theorem~\ref{Ln stabilize}). This integer
is related to the reduction number of a certain ideal. We begin our
exploration by determining the reduction numbers of the ideals
$(J,b)$, where $b$ is a general element in $I$ and $J$ is a
reduction of $I$.

\begin{definition}\label{sdef}{\rm Let $R$ be a Noetherian local ring.
Let $I$ be an $R$--ideal and $J$ a reduction of $I$. Assume $I=(f_1,
\cdots ,f_m)$ and write $T=R[X_{1}, \cdots, X_{m}]$, where $X_{1},
\cdots, X_{m}$ are variables over $R$. Let $\widetilde{K}=(J,{\ds
\sum_{j=1}^{m}}X_{j}f_{j})$. We define the integer $s$ to be
$r_{JT_{\f{m}T}}(\widetilde{K}_{\f{m}T})$.

Notice that since $J$ is a reduction of $I$ it follows that $JT$ is
a reduction of $IT$. Hence $a={\ds \sum_{j=1}^{m}}X_{j}f_{j} \in IT$
is integral over $JT$ and thus $JT$ is a reduction of
$\widetilde{K}$.}
\end{definition}

\s

\begin{lemma} \label{red no leq}
Let $R$ be a Noetherian local ring with infinite residue field. Let
$I$ be an $R$--ideal and $J$ a reduction of $I$. Let $s$ be the
integer as in Definition~{\rm{\ref{sdef}}}. If $b$ is a general
element in $I$, then $r_{J}((J,b)) \leq s$.
\end{lemma}

\proof

Let $k$ denote the residue field of $R$. Let
$\widetilde{M}=\widetilde{K}^{s+1}/J\widetilde{K}^s$,
$\underline{X}=[X_{1},\cdots, X_{m}]$, and
$\underline{\lambda}=[\lambda_{1},\cdots, \lambda_{m}] \in
\mathbb{A}_{R}^{m}$. Write $\f{m}_{\underline{\lambda}}$ for the
maximal ideals $(\f{m},\underline{X}-\underline{\lambda})$ of $T$
and consider the $R$--homomorphisms $\pi_{ \underline{\lambda}} : T
\rightarrow R$ that send $\underline{X}$ to $ \underline{\lambda}$.

From the choice of $s$ we have that $ \widetilde{M}_{\f{m}T}=0$ and
hence $\f{m}T \notin {\rm{Supp}} (\widetilde{M})$. Thus there exists
a dense open subset $U \subset \mathbb{A}_{k}^{m}$ such that
$\widetilde{M}_{\f{m}_{\underline{\lambda}}}=0$ for every $
\ol{\underline{\lambda}} \in U$, where $\ol{\underline{\lambda}}$
denotes the image of ${\underline{\lambda}}$ in
$\mathbb{A}_{k}^{m}$. Therefore for all $\ol{\underline{\lambda}}
\in U$
\begin{equation}\label{(10)}
\widetilde{K}^{s+1}_{\f{m}_{\underline{\lambda}}} =(J
\widetilde{K}^{s})_{\f{m}_{\underline{\lambda}}}.
\end{equation}
In addition we consider the evaluation maps
$\rho_{\underline{\lambda}}: T_{\f{m}_{\underline{\lambda}}}
\rightarrow R$ that send $\underline{X}$ to $ \underline{\lambda}$.
Then for every $\underline{\overline{\lambda}} \in U$ we have
$\rho_{\underline{\lambda}}(\widetilde{K}^{s+1}_{\f{m}_{\underline{\lambda}}})=\rho_{\underline{\lambda}}((J
\widetilde{K}^{s})_{\f{m}_{\underline{\lambda}}})$ according to
equation (\ref{(10)}). In other words $(J,b)^{s+1}=J(J,b)^{s}$,
whenever $b={\ds \sum_{j=1}^{m} } \lambda_{j}f_{j}$,
$\ul{\lambda}=[\lambda_1, \cdots, \lambda_m]$, and
$\ol{\ul{\lambda}} \in U$. Thus $r_{J}((J,b)) \leq s$. \QED

\s

The integer $s$ is in general difficult to compute. However if the
ideal $I$ is $\f{m}$--primary then the following proposition gives a
way to compute this integer.

\begin{proposition} \label{red no eq}
Let $R$ be a Noetherian local ring that is an epimorphic image of a
Cohen--Macaulay ring. Let $\f{m}$ be the maximal ideal of $R$ and
assume that $k=R/\f{m}$ is infinite. Let $I$ be an $\f{m}$--primary
ideal and $J$ a reduction of $I$. Then $r_{J}((J,b)) = s$, where $b$
is a general element in $I$.
\end{proposition}

\proof According to Lemma \ref{red no leq} we have that
$JT_{\f{m}T}$ is a reduction of $\widetilde{K}_{\f{m}T}$ and
$r_{J}((J,b)) \leq s$. Following the notation in the proof of Lemma
\ref{red no leq} we write $\f{m}_{\underline{\lambda}}$ for the
maximal ideals $(\f{m},\underline{X}-\underline{\lambda})$ of $T$,
$\underline{X}=[X_{1}, \cdots, X_{m}]$, and
$\underline{\lambda}=[\lambda_{1}, \cdots, \lambda_{m}] \in
\mathbb{A}_{R}^{m}$. Let $k$ denote the residue field of $R$.

As $J^{s} \subset \widetilde{K}^{s}$ and $J$ is $\f{m}$--primary we
conclude that $(T/\widetilde{K}^{s})_{(\f{m}T)}$ is Artinian and
thus Cohen--Macaulay. By the openness of the Cohen--Macaulay locus
{\cite[Theorem 24.5]{Ma}} there exists a dense open subset $U$ in
$\mathbb{A}_{k}^{m}$ such that
$(T/\widetilde{K}^{s})_{\f{m}_{\underline{\underline{\lambda}}}}$ is
Cohen--Macaulay for all $\ol{\underline{\lambda}} \in U$, where
$\ol{\underline{\lambda}}$ denotes the image of
${\underline{\lambda}}$ in $\mathbb{A}_{k}^{m}$. Write
$b_{\underline{\lambda}}={\ds \sum_{j=1}^{m} } \lambda_jf_{j}$ for
$\ol{\underline{\lambda}} \in U$ and let $W=\{\;b\;|\;
b=b_{\underline{\lambda}}\; \; \mbox{\rm for some} \;\;
\underline{\lambda} \in \mathbb{A}_{R}^{m} \; \; \mbox{\rm such
that} \; \; \ol{\underline{\lambda}} \in U\; \}$.

Suppose that $(J,b)^{s}=J(J,b)^{s-1}$ for some $b \in W$. Notice
that $\sqrt{\widetilde{K}^{s}}=\f{m}T$ and thus ${\rm{dim}} \,
(T/\widetilde{K}^{s})_{\f{m}_{\underline{\lambda}}}=m$. Note that
$\ds{\frac{(T/\widetilde{K}^{s})_{\f{m}_{\underline{\lambda}}}}
{(\underline{X}- \underline{\lambda})}} \simeq R/(J,b)^{s}$ is an
Artinian ring since $I$ is $\f{m}$--primary. In addition
$(T/\widetilde{K}^{s})_{\f{m}_{\underline{\lambda}}}$ is a
Cohen--Macaulay ring of dimension $m$ and the sequence
$\ul{X}-\ul{\lambda}$ consists of $m$ elements in
$\f{m}_{\underline{\lambda}}$. Hence $\underline{X}-
\underline{\lambda}$ is a regular sequence on
$(T/\widetilde{K}^{s})_{\f{m}_{\underline{\lambda}}}$. Consider the
$R$--homomorphisms
$\rho_{\underline{\lambda}}:T_{\f{m}_{\underline{\lambda}}}
 \rightarrow R$ that send $\underline{X}$ to
 $\underline{\lambda}$. Then since $(J,b)^{s}=J(J,b)^{s-1}$ we have
$\rho_{\underline{\lambda}}({\widetilde{K}^{s}}_{\f{m}_{\underline{\lambda}}})
=\rho_{\underline{\lambda}}((J\widetilde{K}^{s-1})_{\f{m}_{\underline{\lambda}}})$.
Thus
\[\begin{array}{cc}
 \widetilde{K}^{s}_{\f{m}_{\underline{\lambda}}}
&=(J\widetilde{K}^{s-1})_{\f{m}_{\underline{\lambda}}}
+(\underline{X}- \underline{\lambda})\cap
\widetilde{K}^{s}_{\f{m}_{\underline{\lambda}}}\\
&=(J\widetilde{K}^{s-1})_{\f{m}_{\underline{\lambda}}}+(\underline{X}-
\underline{\lambda})
\widetilde{K}^{s}_{\f{m}_{\underline{\lambda}}},
\end{array}\]
where the last equality holds since $\underline{X}-
\underline{\lambda}$ is a regular sequence on
$(T/\widetilde{K}^{s})_{\f{m}_{\underline{\lambda}}}$.

Finally by Nakayama's lemma
$\widetilde{K}^{s}_{\f{m}_{\underline{\lambda}}}=(J\widetilde{K}^{s-1})_{\f{m}_{\underline{\lambda}}}$
and therefore
$\widetilde{K}^{s}_{\f{m}T}=(J\widetilde{K}^{s-1})_{\f{m}T}$, which
is a contradiction. \QED

\s

The following theorem makes an effective use of the integer $s$ as
in Definition~\ref{sdef}.

\begin{theorem} \label{Ln stabilize}
Let $R$ be a local Gorenstein ring with infinite residue field $k$.
Let $I$ be an $R$--ideal with $g=\height I
>0$ and $\ell=\ell(I)$, and let $J$ be a minimal reduction
of $I$. Assume that $I$ satisfies $G_{\ell}$ and ${\rm{depth}} \,
R/I^{j} \geq {\rm{dim}} \, R/I-j+1$ for $1\leq j \leq \ell-g$. For
every positive integer $n$ write $L_n=L_{n}(J,I)=
J^{n+1}:K_{n}(J,I)$. Then
\[
L_n=L_s
\]
for every $n \geq s$, where $s$ is as in
Definition~{\rm{\ref{sdef}}}.
\end{theorem}

\proof Let $n$ be a fixed positive integer such that $n \geq s$. By
Theorem~\ref{Kn general} and by Lemma~\ref{red no eq} there exists a
positive integer $t$ such that $K_{n}={\ds \sum_{i=1}^{t}
}(J,b_{i})^n$, where $b_1, \cdots, b_{t}$ are general elements in
$I$, and $r_{J}((J,b_{i})) \leq s$ for all $1 \leq i \leq t$. For
simplicity we denote $(J,b_{i})$ by $J_{i}$ and $K_{n}(J)$ by
$K_{n}$ since $J$ is fixed. Notice that for all $1 \leq i \leq t$ we
have that $J_{i}^{\;s+1}=JJ_{i}^{s}$ according to Lemma~\ref{red no
leq}. Then
\[
K_{n}={\ds \sum_{i=1}^{t} } J_{i}^{n} ={\ds \sum_{i=1}^{t} }
J^{n-s}J_{i}^{s}= J^{n-s} {\ds \sum_{i=1}^{t} } J_{i}^{s} \subset
J^{n-s} K_{s}.
\]In general $J^{n-s}K_{s} \subset K_{n}$ and thus
\begin{equation}\label{(11)}
J^{n-s}K_{s}= K_{n}.
\end{equation}
In conclusion
\[\begin{array}{lcl}
L_{n}=J^{n+1}:K_{n} &\stackrel{(\ref{(11)})}{=}&
J^{n+1}:J^{n-s}K_{s}\\&=& (J^{n+1}:J^{n-s}):K_{s}
\\&\stackrel{(1)}{=}& J^{s+1}:K_{s}=L_{s},\end{array}
\]
where (1) holds since the associated graded ring, ${\rm gr}_{J}(R)$,
of $J$ is Cohen-Macaulay (\cite[Theorem 9.1]{HSV}) and $\height J
>0$ (c.f. \cite[Remark 4.3]{PU}). \QED

\section{Examples}{\label{Examples}}

Finally we arrrive at our goal. We are now ready to answer
Question~\ref{core=Ln?} with the next example using the results from
the previous sections and the computer algebra program
Macaulay~2~(\cite{M2}).

\begin{example} \label{mainexample}
{\rm Let $R=k[x,y,z]_{(x,y,z)}/(z^3)$, where $k$ is an infinite
field of characteristic $2$. Consider the $R$--ideal
$I=(x^2,y^2,xz,yz)$. Then
\begin{enumerate}
\item[{\rm (a)}] $R$ is a 2--dimensional local Gorenstein ring with maximal ideal
$\f{m}=(x,y,z)R$;
\item[{\rm (b)}] $I$ is an $\f{m}$--primary ideal;
\item[{\rm (c)}] $g=\height I=2$, $\ell=\ell(I)=2$, $r=r(I)=2$, and
$r-\ell+g=2$.
\end{enumerate}
We claim that
\[
J^{n+1}:I^{n} \subsetneq {\rm{core}}(I) \subsetneq J^{n+1}:\ds
\sum_{b \in I} (J,b)^n
\]
for any general minimal reduction $J$ of $I$ and any positive
integer $n$.

The computation of $\core{I}$ with Macaulay 2~({\cite{M2}}) is done
using general minimal reductions as in \cite[Theorem~4.5]{CPU01}.
That is, $\core{I} =\ds{\bigcap_{i=1}^{\gamma(I)} J_{i}}$, where
$J_{1}, \cdots, J_{\gamma(I)}$ are general minimal reductions. The
sequence of ideals $\{J^{n+1}:I^{n}\}_{n \in \mathbb{N}}$ is a
decreasing sequence and it stabilizes for $n \geq {\rm{max}}
\{r_J(I)-\ell(I)+g,0\}=2$, according to \cite[Corollary 2.3]{PU}.
Also, $J^{\;n+1}:I^{n} \subset \core{I}$ for $n \geq {\rm{max}}
\{r_J(I)-\ell(I)+g,0\}=2$, according to Theorem~\ref{geninclu}.
Hence it is enough to consider $J^{n+1}:I^{n}$ for $n \leq 2$. Using
Macaulay~2~(\cite{M2}) it is easy to check that $\core{I} \neq
J^{\;n+1}:I^{n}$ for $n \leq 2$, where $J$ is a general minimal
reduction of $I$. Therefore
\[
\core{I} \neq J^{n+1}:I^{n}
\]
for any general minimal reduction $J$ of $I$ and every positive
integer $n$.

Notice that Theorem~\ref{Kn general} provides an algorithm for
computing the ideals $K_{n}$ for any positive integer $n$. Once we
obtain these ideals we can compute $L_{n}(J)=J^{n+1}:K_{n}(J)$. By
Proposition \ref{red no eq} we have that $s=r_{J}((J,b))$, where $b$
is a general element in $I$. In this case $s=2$. By Theorem~\ref{Ln
stabilize} we have that the sequence of the ideals $L_{n}(J)$
stabilizes after $s$ steps. We then check that $\core{I} \neq
L_{n}(J)$ for $n \leq s=2$ and therefore conclude that
\[
\core {I} \neq L_{n}(J)
\]
for all positive integers $n$ and any general minimal reduction $J$
of $I$.

In order to see how close the $\core{I}$ and the ideal $L_{n}(J)$
are we give a description in terms of generating sets obtained using
Macaulay~2~(\cite{M2}). Note that the monomial ideal
$J=(x^{2},y^{2})$ is a minimal reduction of $I$. Then
$\core{I}=(x^2z^2,y^2z^2,x^4,y^4,x^3yz,xy^3z,x^2y^2z,x^2y^3,x^3y^2)$
and $L_{2}(J)=(x^2y^2,y^2z^2,x^4,y^4,x^3yz,xy^3z,x^2y^2)$. Clearly
$x^2y^2 \in L_{2}(J)$ and $x^2y^2 \not \in \core{I}$.}

\end{example}

\s

The question still remains: what is $\core{I}$? According to
\cite[Theorem~4.5]{CPU01} in order to compute the core of an ideal
$I$ we only need to consider a finite intersection of general
minimal reductions. Let $\gamma (I)$ be the number of reductions
required in this intersection. Polini and Ulrich prove that the core
is always contained in the ideals $L_{\;n}(J)$ for every $n$ and any
minimal reduction $J$ of $I$ (\cite[Theorem 4.4]{PU}). On the other
hand $L_{n}(J) \subset J$ for every $n$ and every ideal $J$.

Combining these results we have
\[
\core I \subset \bigcap_{i=1}^{\gamma(I)} L_n(J_i) \subset
\bigcap_{i=1}^{\gamma(I)}J_i =\core I,
\]where $J_1, \cdots, J_{\gamma(I)}$ are general minimal reductions of $I$.
 Therefore
\[
\core I = \bigcap_{i=1}^{\gamma(I)} L_n(J_i),
\]where $J_1, \cdots, J_{\gamma(I)}$ are general minimal reductions of $I$.
In practice though it seems one can do much better. We test this in
Example~\ref{mainexample} using Macaulay~2~(\cite{M2}):

\begin{example}{\rm
In the case of Example \ref{mainexample} we have that for any
positive integer $n \geq s=2$
\[
\core {I} \neq L_{n}(J) \hspace{1cm} {\rm but}  \hspace{1cm}
\core{I} ={\ds \bigcap_{j=1}^2 } L_{n}(J_{j}),
\]
where $J, J_{1},$ and $J_{2}$ are general minimal reductions of $I$,
and $s$ is as in Definition~\ref{sdef}.}
\end{example}

\begin{remark}{\rm
Notice that in the above example, $\ell=2$ and $\gamma(I) \gg 0$. Using Macaualy 2
(\cite{M2}) we were able to construct a series of examples that
support the following conjecture for the core of an ideal in arbitrary characteristic:}
\end{remark}

\begin{conjecture}{\label{conj}}{\rm
Let $R$ be a local Gorenstein ring with infinite residue field, let
$I$ be an $R$--ideal with $g=\height I > 0$ and $\ell=\ell(I)$.
Assume $I$ satisfies $G_{\ell}$ and $\depth R/I^j \geq {\rm dim} \,
R/I -j+1$ for $1 \leq j \leq \ell-g$. Let $s$ be as in
Definition~\ref{sdef}. Then for any integer $n \geq s$
\[
\core{I} ={\ds \bigcap_{j=1}^{\ell}} L_{n}(J_{j}),
\]
where $J_{1}, \cdots, J_{\ell}$ are general minimal reductions of
$I$.}
\end{conjecture}

The above conjecture is consistent with the result of Polini and
Ulrich in the case of analytic spread one (Theorem~\ref{ell=one}).
Notice that Conjecture~\ref{conj} implies that we need to intersect
fewer of the ideals $L_{n}(J)$, since in general $\ell(I) \ll
\gamma(I)$. The following examples support this conjecture.

\begin{example}\label{example3}{\rm
Let $R=k[x,y,z]_{(x,y,z)}/(z^5)$, where $k$ is an infinite field of
characteristic $2$. Consider the $R$--ideal $I=(xy^2z^2+y^5, xyz^3,
x^4y+x^5, x^3yz)$. Then
\begin{enumerate}
\item[{\rm (a)}] $R$ is a 2--dimensional local Gorenstein ring with maximal ideal
$\f{m}=(x,y,z)R$;
\item[{\rm (b)}] $I$ is an $\f{m}$--primary ideal;
\item[{\rm (c)}] $g=\height I=2$, $\ell=\ell(I)=2$, $r=r(I)=4$, and
$r-\ell+g=4$.
\end{enumerate}
Using the same methods as in Example \ref{mainexample} we compute
$\core{I}$ and the ideals $J^{n+1}:I^{n}$ and $L_{n}(J)$ for all $n
\leq s=4$, where $s$ is as in Definition~\ref{sdef}. We conclude
that
\[
J^{n+1}:I^{n} \subsetneq {\rm{core}}(I) \subsetneq L_{n}(J),
\]
for any general minimal reduction $J$ of $I$ and any positive
integer $n$. Nevertheless for all $n \geq 4$
\[
\core{I} ={\ds \bigcap_{j=1}^{2} } L_{n}(J_{j}),
\]
where $J_{1}, J_{2}$ are general minimal reductions of $I$. This is
consistent with Conjecture \ref{conj}.

Since $r-\ell+g=4$ we may repeat the same computations with $k$ an
infinite field of characteristic $3$. Using Macaulay 2 ({\cite{M2}})
we obtain $K_n=I^n$ for $n \geq 4$ and thus
\[
\core{I}= J^{n+1}:I^{n}=L_{n}(J),
\]for any minimal reduction $J$ of $I$ and $n \geq 4$.
Notice that this does not contradict Conjecture~\ref{conj}.}
\end{example}

In both Example \ref{mainexample} and Example \ref{example3} the
analytic spread is 2. We now consider an example where the analytic
spread is 3.

\begin{example}\label{example4}{\rm
Let $R=k[x,y,z,w]_{(x,y,z,w)}/(w^3)$, where $k$ is an infinite field
of characteristic $2$. Consider the $R$--ideal
$I=(x^5,x^2y^2w+z^5,xy^2w^2+x^2z^2w,xyz^2w+y^5,y^2z^2w)$. Then
\begin{enumerate}
\item[{\rm (a)}] $R$ is a 3--dimensional local Gorenstein ring with maximal ideal
$\f{m}=(x,y,z,w)R$;
\item[{\rm (b)}] $I$ is an $\f{m}$--primary ideal;
\item[{\rm (c)}] $g=\height I=3$, $\ell=\ell(I)=3$, $r=r(I)=2$, and
$r-\ell+g=2$.
\end{enumerate}
We again use the same methods as in Example \ref{mainexample} to
compute $\core{I}$ and the ideals $J^{n+1}:I^{n}$ and $L_{n}(J)$ for
all $n \leq s=2$, where $s$ is as in Definition~\ref{sdef}. We
conclude that
\[
J^{n+1}:I^{n} \subsetneq {\rm{core}}(I) \subsetneq L_{n}(J)
\]
for any general minimal reduction $J$ of $I$ and any positive
integer $n$. Nevertheless for all $n \geq 2$
\[
\core {I} \neq {\ds \bigcap_{i=1}^{2}}L_{n}(J_{i}) \hspace{1cm}
{\rm{but}} \hspace{1cm}  \core {I}={\ds \bigcap_{i=1}^{3}}
L_{n}(J_{i}),
\]
where $J_{1}, J_{2}, J_{3}$ are general minimal reductions of $I$.
Thus this example provides yet more evidence for the truth of
Conjecture~\ref{conj}.}
\end{example}

\s

In the case of analytic spread 4 we exhibit the following example.
The computations become quite difficult for higher analytic spreads.

\s

\begin{example}{\label{example5}}{\rm

Let $R=k[x,y,z,w,t]_{(x,y,z,w,t)}/(w^3)$, where $k$ is an infinite
field of characteristic $2$. Consider the $R$--ideal
$I=(x^5,x^2y^2w+z^5,xt^2w^2+x^2z^2w,xyt^2w+y^5,yz^2wt,t^5)$. Then
\begin{enumerate}
\item[{\rm (a)}] $R$ is a 4--dimensional local Gorenstein ring with maximal ideal
$\f{m}=(x,y,z,w,t)R$;
\item[{\rm (b)}] $I$ is an $\f{m}$--primary ideal;
\item[{\rm (c)}] $g=\height I=4$, $\ell=\ell(I)=4$, $r=r(I)=2$, and
$r-\ell+g=2$.
\end{enumerate}
Once again we use the same methods as in Example \ref{mainexample} to
compute $\core{I}$ and the ideals $J^{n+1}:I^{n}$ and $L_{n}(J)$ for
all $n \leq s=2$, where $s$ is as in Definition~\ref{sdef}. We
conclude that
\[
J^{n+1}:I^{n} \subsetneq {\rm{core}}(I) \subsetneq L_{n}(J)
\]
for any general minimal reduction $J$ of $I$ and any positive
integer $n$. Nevertheless for all $n \geq 2$
\[
\core {I} \neq {\ds \bigcap_{i=1}^{2}}L_{n}(J_{i}), \hspace{1cm} \core {I} \neq {\ds \bigcap_{i=1}^{3}}L_{n}(J_{i}), \hspace{1cm} {\rm{but}} \hspace{1cm}  \core {I}={\ds \bigcap_{i=1}^{4}}
L_{n}(J_{i}),
\]
where $J_{1}, J_{2}, J_{3}, J_{4}$ are general minimal reductions of
$I$. Again the validity of  Conjecture~\ref{conj} is supported by
this example.}
\end{example}
\s

In all the previous examples the rings that we considered were non--reduced. The next example is set in a regular local ring.

\begin{example}{\label{sharp}}(\cite{FPU}){\rm \; Let $R=k[x,y]_{(x,y)}$, where $k$ is an infinite perfect field. Let
$I=(x^2y^8,y^9,x^5y^3,x^4y^4,x^6)$. Then
\begin{enumerate}
\item[{\rm (a)}] $R$ is a $2$--dimensional regular local ring with maximal ideal
$\f{m}=(x,y)R$;
\item[{\rm (b)}] $I$ is an $\f{m}$--primary ideal;
\item[{\rm (c)}] $g=\height I =2$, $\ell=\ell(I)=2$, $r=r(I)=2$,
and $r-\ell+g=2$;
\end{enumerate}

Once again following the same ideas as before we conclude that for all $n \geq s=2$
\[
\core{I} \neq L_{n}(J) \hspace{1cm} {\mbox{ \rm{and}}} \hspace{1cm} \core{I}={\ds \bigcap_{j=1}^{2} } L_{n}(J_{j}),
\]
where $J$, $J_{1}$, and $J_{2}$ are general minimal reductions of $I$ and $s$ is as in Definition~\ref{sdef}.}
\end{example}

\s

\section{Acknowledgments}
This paper is part of my Ph.D. thesis which was written under the
supervision of Professor Bernd Ulrich at Purdue University. I would
like to express my gratitude to Professor Bernd Ulrich for his
guidance, advice, and constant support. I would also like to thank
Professor Raymond C. Heitmann (University of Texas, Austin) for
reading the manuscript and for his suggestions that improved the
exposition of the paper.

\s


\begin{thebibliography}{99}
\bibitem{AvHer}{L. Avramov and J. Herzog, The Koszul algebra of a codimension $2$ embedding,  Math. Z.  175  (1980), 249--260.}

\bibitem{CPU01}{A. Corso, C. Polini and B. Ulrich, The structure
of the core of ideals, Math. Ann. 321 (2001), 89--105.}

\bibitem{CPU02}{A. Corso, C. Polini and B. Ulrich, Core and residual intersections of ideals, Trans. Amer. Math. Soc.
354 (2002), 2579--2594.}


\bibitem{M2}{D. Grayson, M. Stillman, Macaulay 2, A computer algebra system for computing in Algebraic Geometry and
Commutative Algebra, available through anonymous ftp from
http://www.math.uiuc.edu/Macaulay2 .}


\bibitem{FPU}{L. Fouli, C. Polini and B. Ulrich, The core of ideals in arbitrary characteristic, preprint.}



\bibitem{HSV}{J. Herzog, A. Simis, and W. V. Vasconcelos, Koszul homology and
blowing--up rings, in Commutative Algebra, Proceedings: Trento 1981
(Greco/Valla eds.), Lecture Notes in Pure and Applied Mathematics
84, Marcel Dekker, New York 1983, 79--169.}

\bibitem{Hun}{C. Huneke, Linkage and the Koszul homology of ideals,  Amer. J. Math.  104  (1982), 1043--1062.}

\bibitem{HS}{C. Huneke and I. Swanson, Cores of ideals in
$2$--dimensional regular rings, Michigan Math. J. 42 (1995),
193--208.}


\bibitem{HT}{C. Huneke and N. Trung, On the core of ideals, Compos. Math.  141  (2005), 1--18.}



\bibitem{HySm1}{E. Hyry and K. Smith, Core versus graded core, and
global sections of line bundles, Trans. Amer. Math. Soc. 356 (2003),
3143--3166.}

\bibitem{HySm2}{ E. Hyry and K. Smith, On a non--vanishing
conjecture of Kawamata and the core of an ideal, Amer. J. Math. 125
(2003), 1349--1410.}



\bibitem{L}{J. Lipman, Adjoints of ideals in regular local rings,
Math. Research Letters 1 (1994), 739--755.}

\bibitem{Ma}{ H. Matsumura, Commutative Ring Theory,  Cambridge Studies in Advanced Mathematics
8, Cambridge University Press, Cambridge, 1989.}

\bibitem{NR}{D.G. Northcott and D. Rees, Reductions of ideals in
local rings, Proc. Camb. Phil. Soc. 50 (1954), 145--158.}



\bibitem{PU}{C. Polini and B. Ulrich, A formula for the core of an ideal,
 Math. Ann.  331  (2005), 487--503.}

\bibitem{PUV}{C. Polini, B. Ulrich and M. Vitulli, The core of zero-dimensional monomial
ideals, Adv. Math. 72 (2007), 72-93.}


\bibitem{RS}{D. Rees and J.D. Sally, General elements and joint
reductions, Michigan Math. J. 35 (1988), 241--254.}


\end{thebibliography}
\end{document}